\documentclass[11pt,reqno,a4paper,english]{amsart}

\usepackage[utf8]{inputenc}  
\usepackage{dsfont}
\usepackage{xcolor}

\usepackage{enumerate}          
\usepackage{amssymb,amsmath,amsfonts,amsthm}            
\usepackage{pdfsync}
\usepackage{url}
\usepackage{hyperref}
\usepackage{soul}
\usepackage{anysize}
\usepackage{mathrsfs}
\usepackage{graphicx,float,color}
\usepackage{epsfig}
\usepackage{tikz}
\usepackage{caption}
\usetikzlibrary{patterns, patterns.meta}
\usetikzlibrary{calc,arrows.meta,positioning}
\usetikzlibrary{decorations.pathreplacing}

\begin{document}
\theoremstyle{plain}
\newtheorem{theorem}{Theorem}[section]
\newtheorem{definition}[theorem]{Definition}
\newtheorem{proposition}[theorem]{Proposition}
\newtheorem{lemma}[theorem]{Lemma}

\newtheorem{corollary}[theorem]{Corollary}
\newtheorem{conjecture}[theorem]{Conjecture}
\newtheorem{remark}[theorem]{Remark}
\newtheorem{open}[theorem]{Open question}
\numberwithin{equation}{section}
\newtheorem*{claim}{Claim}
\newtheorem{assumption}[theorem]{Assumption}
\errorcontextlines=0

\newcommand{\Id}{\text{Id}}
\newcommand{\mnum}{A}
\newcommand{\mrat}{B}
\newcommand{\R}{\mathbb{R}}
\newcommand{\Z}{\mathbb{Z}}
\newcommand{\N}{\mathbb{N}}
\newcommand{\X}{\mathcal{X}}
\newcommand{\Y}{\mathcal{Y}}
\newcommand{\supp}{\text{supp}}
\newcommand{\red}{\color{red}}
\newcommand{\black}{\color{black}}
\newcommand{\A}{A}
\newcommand{\ppow}{\beta}
\newcommand{\gau}{\ppow}
\newcommand{\Cor}{\mathcal{C}}
\renewcommand{\H}{\mathcal{H}}
\newcommand{\B}{B}
\newcommand{\Ball}{\mathcal{B}}
\newcommand{\C}{C}
\newcommand{\D}{D}
\newcommand{\E}{E}
\newcommand{\ct}{5\B\sqrt{d}}
\newcommand{\pow}{p}
\newcommand{\powun}{p_1}
\newcommand{\powde}{p_2}
\newcommand{\pds}{\delta}
\renewcommand{\P}{\mathbb{P}}
\renewcommand{\S}{\mathcal{S}}
\newcommand{\dist}{{\rm dist}}
\newcommand{\spt}{{\rm spt}}
\newcommand{\diam}{{\rm diam}}
\newcommand{\Var}{{\rm Var}}
\newcommand{\inte}{{\rm int}}
\newcommand{\Per}{{\rm Per}}
\newcommand{\vol}{{\rm vol}}
\newcommand{\CDM}{C_{\rm DM}}
\newcommand{\Res}{\mathcal{B}}
\newcommand{\Kant}{\mathcal K}
\newcommand{\YSp}{\mathcal Y}
\newcommand{\Rsp}{\mathbb R}
\newcommand{\eps}{\varepsilon}
\newcommand{\Class}{\mathcal C}
\newcommand{\Conv}{\mathrm{Conv}}
\newcommand{\sca}[2]{\langle #1\vert #2\rangle}
\newcommand{\nr}[1]{\Vert #1\Vert}
\newcommand{\restr}[1]{\left. #1\right\vert}
\newcommand{\dd}{{\rm d}}

\hyphenation{mar-gi-nals}

\setcounter{tocdepth}{1}

\title{Unstable optimal transport maps}
\author{Cyril Letrouit}\address{Cyril Letrouit. Université Paris-Saclay, CNRS, Laboratoire de mathématiques d’Orsay, 91405, Orsay, France} \email{cyril.letrouit@universite-paris-saclay.fr}
\date{\today}
\maketitle

\begin{abstract}
The stability of optimal transport maps with respect to perturbations of the marginals is a question of interest for several reasons, ranging from the justification of the linearized optimal transport framework to numerical analysis and statistics. Under various assumptions on the source measure, it is known that optimal transport maps are stable with respect to variations of the target measure. 

In this note, we focus on the mechanisms that can, on the contrary, lead to \emph{instability}. We identify two of them, which we illustrate through examples of absolutely continuous source measures $\rho$ in $\R^d$ for which optimal transport maps are less stable, or even very unstable. We first show that instability may arise from the unboundedness of the density: we exhibit a source density on the unit ball of $\R^d$ which blows up superpolynomially at two points of the boundary and for which optimal transport maps are highly unstable. Then we prove that even for uniform densities on bounded open sets, optimal transport maps can be rather unstable close enough to configurations where uniqueness of optimal plans is lost.
%if the support is too disconnected.
\end{abstract}

\section{Introduction}

\subsection{Main results.} Let $\mathcal{P}_2(\R^d)$ denote the set of probability measures on $\R^d$ with finite second moment. Given $\rho,\mu\in\mathcal{P}_2(\R^d)$ such that $\rho$ is absolutely continuous with respect to the Lebesgue measure, Brenier's theorem (\cite{brenier}, \cite{brenierCPAM}) guarantees the existence and $\rho$-a.e. uniqueness of an optimal transport map $T_\mu\in L^2(\rho)$ from $\rho$ to $\mu$. Brenier \cite{brenierCPAM} also proved that the map $\mu\mapsto T_\mu$ is continuous from $(\mathcal{P}_2(\R^d),W_2)$ to $L^2(\rho)$, where $W_p$ denotes the $p$-Wasserstein distance. This result can be interpreted as a statement of (qualitative) stability with respect to perturbations of the target measure, for a fixed source measure $\rho$.

In recent years, several people have tried to quantify this stability, i.e., to provide conditions on $\rho$ that guarantee quantitative bounds on the $L^2$-distance between $T_\mu$ and $T_\nu$ in terms of Wasserstein distances between $\mu$ and $\nu$. Under various assumptions on $\rho$, quantitative stability bounds  have been established (see the literature review in Section \ref{s:litreview}). These are bounds of the form
\begin{equation}\label{e:stabmapoule}
\forall \mu,\nu\in\mathcal{P}_2(\Y), \qquad \|T_\mu-T_\nu\|_{L^2(\rho)} \leq CW_p(\mu,\nu)^{\alpha}
\end{equation} 
 for $C,\alpha>0$, $p\geq 1$ and $\Y\subset\R^d$. There are multiple motivations for looking at this type of inequalities, notably the justification of the linearized optimal transport framework (see below) and the statistical estimation of optimal transport maps (\cite[Chapter 3.2]{chewi}, \cite{balakrishnan}).

However, it has never been proven that optimal transport maps could be \emph{unstable}. In this note we first construct an absolutely continuous source measure $\rho$ for which optimal transport maps are \emph{highly unstable}. More precisely, in Theorem \ref{t:highlyunstable}, we construct $\rho$ whose density is bounded below over the (closed) unit ball $B_{\R^d}(0,1)$ such that for any ball $\Y=B_{\R^d}(0,R)$, any $C,\alpha>0$ and $p\geq 1$, the bound \eqref{e:stabmapoule} fails. In the sequel, $\mathcal{P}(\X)$ denotes the set of probability measures on $\X\subset\R^d$.

%The small difference between these two statements is that in the second one, we want to disprove \eqref{e:stabmapoule} using only measures $\mu,\nu$ in a compact set, and thus we do not allow perturbations ``at infinity". Therefore, Theorem \ref{t:highlyunstable} is in a sense more difficult to obtain, and our construction requires a source measure $\rho$ whose density blows up at the boundary of its (compact) support -- whereas in Theorem \ref{t:noncompacttarget} we may take $\rho$ uniform, and in particular bounded above and below. 

%\begin{theorem}[$\rho$ uniform, targets over $\R^d$]\label{t:noncompacttarget}Let $d\geq 2$. There exists a bounded open subset $\X\subset\R^d$ such that if $\rho$ denotes the uniform probability measure on $\X$, then for any $C,\alpha>0$ and $p\geq 1$, the inequality\begin{equation}\label{e:stabmap}\forall \mu,\nu\in\mathcal{P}_2(\R^d), \qquad \|T_\mu-T_\nu\|_{L^2(\rho)} \leq CW_p(\mu,\nu)^{\alpha}.\end{equation} fails. \end{theorem}

\begin{theorem}\label{t:highlyunstable}
Let $d\geq 2$. There exists  an absolutely continuous $\rho\in\mathcal{P}(B_{\R^d}(0,1))$ with density bounded below, such that for any ball $\Y=B_{\R^d}(0,R)$ with $R>0$, any $C,\alpha>0$ and $p\geq 1$, the inequality
\begin{equation}\label{e:stabmap}
\forall \mu,\nu\in\mathcal{P}(\Y), \qquad \|T_\mu-T_\nu\|_{L^2(\rho)} \leq CW_p(\mu,\nu)^{\alpha}
\end{equation} 
fails.
\end{theorem} 

The density $\rho$ we construct to prove Theorem \ref{t:highlyunstable} blows up slightly superpolynomially at two points of the boundary of $B_{\R^d}(0,1)$.  The results of \cite{letrouitmerigot} show that if $\rho$ blows up only polynomially at the boundary of $B_{\R^d}(0,1)$, then \eqref{e:stabmap} holds for some $C,\alpha>0$ and $p=1$ (see Remark \ref{r:superpolyn}).

Apart from unbounded densities, if we focus on source densities bounded above and below on bounded open sets, some loss of stability can nevertheless occur close enough to configurations where optimal plans are non-unique. It is actually natural to expect loss of stability just before loss of uniqueness. To illustrate this idea, we show the following result, using as a source measure the \emph{uniform probability measure} $\rho$ on a carefully chosen \emph{bounded open set} $\X\subset\R^d$:
\begin{theorem}\label{t:unstableW1}
Let $d\geq 2$. There exists a bounded open set $\X\subset\R^d$ and a compact set $\Y$ such that if $\rho$ denotes the uniform probability measure on $\X$, then for any $C>0$, $p\geq 1$ and any $\alpha>\frac{p}{2(p+1)}$ (in particular for $\alpha=\frac12$), the inequality
\begin{equation}\label{e:stabmap}
\forall \mu,\nu\in\mathcal{P}(\Y), \qquad \|T_\mu-T_\nu\|_{L^2(\rho)} \leq CW_p(\mu,\nu)^{\alpha}.
\end{equation} 
fails. 
\end{theorem}
We do not know if it would be possible to obtain the (stronger) conclusion of Theorem \ref{t:highlyunstable} for some $\rho$ bounded above and below on a well-chosen bounded open set $\X$. Necessarily, $\X$ would need to have an infinite number of connected components, see Remark \ref{r:finiteunion}.

Our results shed light on the regularity of the map $\mu\mapsto T_\mu$ from $(\mathcal{P}(\Y),W_2)$ to $L^2(\rho)$, which reflects the stability of optimal transport maps with source measure $\rho$ with respect to perturbations of the marginals.  Gigli \cite[Theorem 5.1]{gigli} proved in 2011 that in some situations, this map is not better than $\frac12$-Hölder. Since then, it remains an open question to determine when this regularity is actually achieved: for instance, Gigli \cite[Corollary 3.4]{gigli} (see also \cite[Theorem 2.3]{merigot}) showed that it is achieved at any $\mu$ such that $T_\mu$ is Lipschitz. This result plays a key role in statistical optimal transport (\cite[Chapter 3.2]{chewi}, \cite{balakrishnan}). No result has shown that $\mu\mapsto T_\mu$ could fail to be $\frac12$-Hölder. In the proof of Theorem \ref{t:unstableW1} (with $p=2$), we show that for $\rho$ uniform on a well-chosen bounded open set, this map is in fact not better than $\frac13$-Hölder at some (explicit) $\mu$. And Theorem \ref{t:highlyunstable} shows that if $\rho$ is allowed to have unbounded density, then it can happen that for any $\alpha>0$ the map is not $\alpha$-Hölder.

% that if $\Y$ is compact and $\rho$ is bounded above and below on a compact set, \eqref{e:stabmapoule} should always hold for $p=2$, $\alpha=\frac12$ and some $C>0$.

\subsection{Previous works.}\label{s:litreview}
Inequalities of the form \eqref{e:stabmap} have been established under various assumptions on $\rho$. After works by Gigli \cite{gigli}, Berman \cite{berman}, Mérigot-Delalande \& Chazal \cite{merigot} and Delalande \& Mérigot \cite{delmer}, the following result has been achieved in \cite{letrouitmerigot}:
\begin{theorem}\cite[Theorem 1.7]{letrouitmerigot} \label{t:avecquentin}
Let $\X\subset\R^d$ be a John domain with rectifiable boundary, and let $\rho$ be a probability density on $\X$, bounded from above and below by positive constants. Then, for any compact set $\mathcal{Y}$, there exists $C_{\rho,\mathcal{Y}}>0$ such that for any probability measures $\mu,\nu$ supported in $\mathcal{Y}$,
\begin{equation}\label{e:stabmapjohn2}
\|T_\mu-T_\nu\|_{L^2(\rho)}\leq C_{\rho,\mathcal{Y}}W_1(\mu,\nu)^{1/6}.
\end{equation}
\end{theorem}
Recall that any bounded and connected Lipschitz domain is a John domain, and therefore Theorem \ref{t:avecquentin} applies for instance in this case. The paper \cite{letrouitmerigot} also establishes similar stability inequalities for log-concave $\rho$, and for $\rho$ blowing-up (or decaying) polynomially at the boundary of a smooth compact set. 
%Following the strategy of \cite{delmer}, the inequality \eqref{e:stabmapjohn2} is established in two steps: first proving an analogous stability inequality for Kantorovich potentials (recall that the Brenier map $T_\mu$ is the gradient of a convex function $\phi_\mu$, called Kantorovich potential), and then deducing the stability of maps using arguments from geometric measure theory. 
Quantitative stability inequalities have also been proved for optimal transport maps with respect to $p$-costs in $\R^d$ \cite{mischler}, and the squared distance cost on Riemannian manifolds \cite{kitagawa}.

Quantitative stability results were primarily motivated by numerical analysis questions: if $\mu$ is known only through an approximation $\nu=\widehat{\mu}$ (for instance through samples), is is true that the optimal transport map $T_{\widehat{\mu}}$, which one may compute for instance through semi-discrete optimal transport, is not far from $T_\mu$? 

Quantitative stability inequalities, when they hold, also serve as a justification for the linearized optimal transport framework introduced in \cite{slepcev} (and used later in several applications): the mapping $\mu \mapsto T_{\mu}$ provides an embedding of $(\mathcal{P}(\Y),W_2)$ into the Hilbert space $L^2(\rho,\R^d)$, and this embedding allows one to apply the standard ``Hibertian" statistical toolbox to measure-valued data. This embedding is distance-increasing, meaning that $\|T_\mu - T_\nu\|_{L^2(\rho)}\geq W_2(\mu,\nu)$, and stability estimates such as \eqref{e:stabmapjohn2} (combined with the fact that $W_1\leq W_2$) show that it is bi-Hölder continuous when $\rho$ satisfies the assumptions of Theorem \ref{t:avecquentin}. In other words, the distance $d(\mu,\nu) = \|T_\mu - T_\nu\|_{L^2(\rho)}$ preserves in a rough way the geometry associated to the Wasserstein distance. However, recall that some results show the impossibility of embedding (in a very coarse sense, in particular in a bi-Hölder way) Wasserstein spaces over $\R^d$ ($d\geq 3$) into Banach spaces of non-trivial type such as Hilbert spaces (see e.g. \cite{andoni}); this is why the literature about quantitative stability inequality has mostly focused on the case where targets are taken over a \emph{compact set} $\Y$.

On the side of instability, the only known results hold merely for \emph{Kantorovich potentials}. Recall that given $\rho,\mu\in\mathcal{P}_2(\R^d)$ with $\rho$ absolutely continuous, a Kantorovich potential is a convex function whose gradient is equal to the Brenier map $T_\mu$. If the support of $\rho$ is connected, then there exists for each $\mu\in\mathcal{P}_2(\R^d)$ a unique Kantorovich potential $\phi_\mu$ satisfying $\int \phi_\mu \dd\rho=0$. The paper \cite{letrouitmerigot} provides examples of uniform source measures $\rho$ on non-John domains $\X$ (but still bounded and connected) for which Kantorovich potentials are highly unstable in the same sense as above: the inequality 
$$
\forall \mu,\nu\in\mathcal{P}(\Y), \qquad \|\phi_\mu-\phi_\nu\|_{L^2(\rho)} \leq CW_p(\mu,\nu)^{\alpha}
$$
fails for any $C,\alpha>0$ and $p\geq 1$. However it does not provide examples of unstable transport maps.

\subsection{Acknowledgments.} The author acknowledges the support of the Agence nationale de la recherche, through the PEPR PDE-AI project (ANR-23-PEIA-0004).

\section{Proof of Theorem \ref{t:highlyunstable}}
Let $A=(1,0,\ldots,0)\in\R^d$, $A'=(-1,0,\ldots,0)\in\R^d$, and $\mathcal{E}=\{A,A'\}$, and let
$$
f:r\mapsto r^{-d}\min(1,(\log r)^{-2})
$$
for $r>0$. The function $g:x\mapsto f(\dist(x,\mathcal{E}))$ is integrable on $B_{\R^d}(0,1)$, where $\dist(x,\mathcal{E})$ denotes the Euclidean distance from $x$ to the set $\mathcal{E}$. To see this, we only need to check integrability close to $\mathcal{E}$. Taking polar coordinates around $A$, $\varepsilon>0$ small, and denoting by $\sigma_{d-1}$ the area of the $(d-1)$-dimensional unit sphere, we have
$$
\int_{B_{\R^d}(A,\varepsilon)\cap B_{\R^d}(0,1)} g(x) \dd x\leq \sigma_{d-1} \int_0^\varepsilon r^{-d}(\log r)^{-2} r^{d-1} \dd r = \sigma_{d-1}\int_0^\varepsilon \frac{1}{r (\log r)^2} \dd r <+\infty.
$$
The same computation holds replacing $A$ by $A'$. Therefore we may choose $c_0>0$ in a way that the density
$$
\rho(x)=c_0f(\dist(x,\mathcal{E}))
$$
on $B_{\R^d}(0,1)$ is a probability density (here and in the sequel, absolutely continuous measures with respect to the Lebesgue measure are identified with thei density).

Let us now consider $B_\theta=(R\sin(\theta),R\cos(\theta),0,\ldots,0)$ and $B_\theta'=(-R\sin(\theta),-R\cos(\theta),0,\ldots,0)$ for $\theta\in\R$. We set 
$$
\mu_\theta=\frac12(\delta_{B_\theta}+\delta_{B_\theta'}).
$$
Since $\rho$ is invariant under the transformation $x\mapsto -x$, the optimal transport map from $\rho$ to $\mu_\theta$ is the map
$$
T_{\mu_\theta}:x\mapsto \begin{cases} B_\theta & \text{  if } \langle x,B_\theta\rangle > 0  \\  B_\theta' & \text{  if } \langle x,B_\theta\rangle < 0.  \end{cases}
$$
In other words $T_{\mu_\theta}$ sends each point of the source to the closest point in the support of the target.

For small $\theta$, most points that $T_{\mu_0}$ sends to $B_0$ are sent to $B_\theta$ under $T_{\mu_\theta}$, but some of them, those satisfying $\langle x,B_0\rangle>0>\langle x,B_\theta\rangle$, are sent to $B_\theta'$. We will bound from below the measure of the latter points in order to bound from below $\|T_{\mu_\theta}-T_{\mu_0}\|_{L^2(\rho)}$. Let us show that for any $\theta\geq 0$ small enough (so that $\sin \theta \geq \theta/2$), 
\begin{equation}\label{e:B0Btheta}
B(A',\theta/4) \cap \{x\in\R^d \mid x_2>0\}\subset \{x\in\R^d \mid \langle x,B_0\rangle>0>\langle x,B_\theta\rangle\}.
\end{equation}
Let $x=(x_1,\ldots,x_d)\in\R^d$ be an element of the left-hand side. Then $\langle x,B_0\rangle>0$ since $x_2>0$. Moreover,
$$
R^{-1}\langle x,B_\theta\rangle = x_1\sin \theta + x_2\cos\theta \leq x_1 \sin \theta + x_2 \leq \Bigl(-1+\frac{\theta}{4}\Bigr)\frac{\theta}{2}+\frac{\theta}{4}<0
$$
which concludes the proof of \eqref{e:B0Btheta}. Let us now observe that there exists $c_1>0$ (independent of $\theta$) such that the left-hand side of \eqref{e:B0Btheta} has measure at least $c_1 \rho(B(A',\theta/4))$. Indeed, for $\theta$ small enough the support of $\rho$ contains the spherical sector
\begin{equation}\label{e:cone}
B(A',\theta/4) \cap \mathcal{C}
\end{equation}
(see Figure \ref{fig:blowup}) where $\mathcal{C}$ is the cone
$$
\mathcal{C}=\left\{x\in\R^d \mid x_2>0, \ \frac{x_1+1}{|x-A'|}>\frac12 \right\}
$$
Denote by $c_1>0$ the angular aperture of this cone relatively to the full solid angle of the $(d-1)$-dimensional unit sphere, i.e., the area of the intersection of $\mathcal{C}$ with the unit sphere $\partial B(A',1)$, divided by the area of the full unit sphere $\partial B(A',1)$.
Since $\rho$ has a radial density close to the center $A'$, the $\rho$-measure of $B(A',\theta/4) \cap \mathcal{C}$ is at least equal to $c_1\rho(B(A',\theta/4))$.
%i.e., the intersection of the left-hand side of \eqref{e:B0Btheta} with a cone with apex at $A'$ and fixed positive aperture.

\begin{center}
\begin{figure}[h!]
\begin{tikzpicture}[scale=2]

    % Axes
    \draw[->] (-2,0) -- (2,0) node[right] {$x_1$};
    \draw[->] (0,-1.5) -- (0,1.5) node[above] {$x_2$};

    % Unit circle
    \draw[thick,blue] (0,0) circle (1);
   % \node at (0.7,0.7) [blue] {$x^2+y^2=1$};

    % Apex of the cone
    \filldraw[orange] (-1,0) circle (0.03) node[below left] {$A'$};

    % Cone lines (tangent to the unit circle)
    % Tangent points: (1/2, ±√3/2)
    \draw[thick,orange] (-1,0) -- (0.5,0.866);  % upper tangent
    \draw[thick,orange] (-1,0) -- (1,0); % lower tangent

    \draw[thick,purple] (-1,0) circle (0.4);

\fill[pattern=north east lines, pattern color=gray]
(-1,0) -- ++(0:0.4) arc[start angle=0, end angle=30, radius=0.4] -- cycle;

    % Labels
   % \node[orange] at (0.2,1) {Cone};
\end{tikzpicture}
\caption{Illustration of \eqref{e:cone}: in purple the ball $B(A',\theta/4)$, in orange the cone in \eqref{e:cone}, in grey dashed lines the set \eqref{e:cone}, and in blue the boundary of the support of $\rho$.} \label{fig:blowup}
\end{figure}
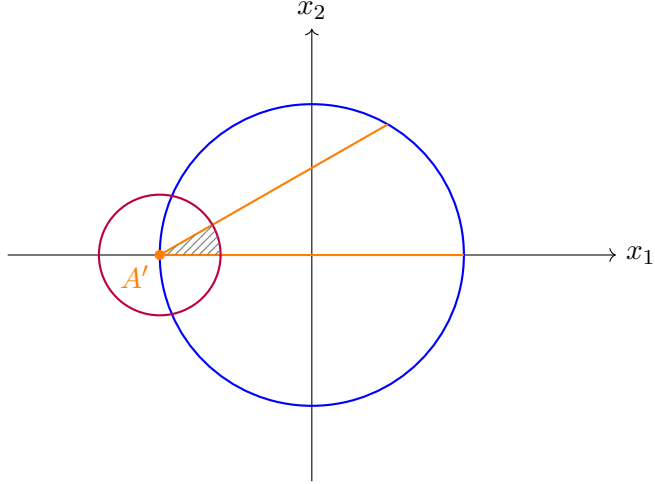
\end{center}

We deduce that for $\theta$ small enough, since $|B'_\theta-B_0|\geq R$,
\begin{align*}
\|T_{\mu_\theta}-T_{\mu_0}\|_{L^2(\rho)}^2 &\geq |B'_\theta-B_0|^2 \rho(\{x \mid \langle x,B_0\rangle>0>\langle x,B_\theta\rangle\})\\
&\geq c_1R^2\rho(B(A',\theta/4))\\
&\geq c_0c_1\sigma_{d-1}R^2\int_0^{\theta/4} \frac{1}{r(\log r)^2} \dd r\\
&=c_0c_1\sigma_{d-1}R^2\frac{1}{|\log(\theta/4)|}
\end{align*}
which decays to $0$ as $\theta\rightarrow 0$ slower than any $\theta^\alpha$ with $\alpha>0$.
Moreover, for $\theta$ small enough, $W_p(\mu_0,\mu_\theta)=R\sin(\theta/2)\sim R\theta/2$. Put together, these estimates conclude the proof of Theorem \ref{t:highlyunstable}.

\begin{remark}\label{r:superpolyn}
The density $\rho$ used in the proof of Theorem \ref{t:highlyunstable} blows up slightly superpolynomially. Let us instead consider $\rho(x)=c_\delta h(\dist(x,\mathcal{E}))$ where 
$$
h:r\mapsto r^{-d+\delta}
$$
for some $\delta>0$, i.e., a density which blows up polynomially at $A$ and $A'$. The assumption $\delta>0$ is necessary to ensure that $\rho$ is integrable, and $c_\delta$ is chosen to make $\rho$ a probability measure. Following the proof strategy of \cite{letrouitmerigot}, notably the proof of Theorem 1.10 in this paper, one can show that for any compact $\Y\subset\R^d$, there exists $C>0$ such that
$$
\|T_\mu-T_\nu\|_{L^2(\rho)}\leq CW_2(\mu,\nu)^{\frac{\delta}{6d}}.
$$
Conversely, the same argument as in the proof of Theorem \ref{t:highlyunstable} (with $h$ replacing $f$) shows that if
\begin{equation}\label{e:polyn}
\|T_\mu-T_\nu\|_{L^2(\rho)}\leq CW_p(\mu,\nu)^{\alpha}
\end{equation}
holds for some $C,\alpha>0$ and $p\geq 1$, then necessarily $\alpha\leq \delta$. 

It is therefore no surprise that for densities which blow up faster than $r^{-d+\delta}$ for any $\delta>0$, any quantitative stability inequality of the form \eqref{e:polyn} breaks down. For these densities, the techniques of \cite{letrouitmerigot} do not apply anymore, and we are not aware of any quantitative stability inequality, even with a larger right-hand side than $W_p(\mu,\nu)^\alpha$ (e.g., $(1+|\log W_p(\mu,\nu)|)^{-1}$).
\end{remark}

\begin{remark}\label{r:finiteunion}
Let us recall that for $\rho$ bounded above and below on a finite union of John domains with rectifiable boundaries, the quantitative stability inequality \eqref{e:stabmapoule} holds for any $p\geq 1$ and $\Y$ compact, with exponent $\alpha=1/6$, as proved in \cite[Remark 4.2]{letrouitmerigot}. Therefore, the conclusion of Theorem \ref{t:highlyunstable} cannot be obtained for $\rho$ bounded above and below, unless the support of $\rho$ has infinitely many connected components. 
\end{remark}

\begin{remark}
For the choice of $\rho$ used in the proof of Theorem \ref{t:highlyunstable}, the strong instability of optimal transport maps naturally raises the following question: how can one statistically estimate optimal transport maps with this source measure?
\end{remark}

\section{Proof of Theorem \ref{t:unstableW1}}

\subsection{Idea of the construction.} We build upon the idea that stability of optimal transport maps can be lost just before uniqueness is lost.  To explain our construction of $\rho$, let us start with a  simple and well-known example where optimal transport plans are not unique: consider the vertices $A,B,A',B'$ of a square, for instance $A=(1,0)$, $B=(0,1)$, $A'=(-1,0)$ and $B'=(0,-1)$, and consider the optimal transport problem from $\frac12 (\delta_A+\delta_{A'})$ to $\frac12(\delta_B+\delta_{B'})$.

\begin{center}
\begin{figure}[h!]
\begin{tikzpicture}[>=stealth,thick]

% Define coordinates
\coordinate (A) at (1,0);
\coordinate (B) at (0,1);
\coordinate (A') at (-1,0);
\coordinate (B') at (0,-1);

% Draw curved arrows between some vertices
\draw[->] (0.9,0.1) to (0.1,0.9);
\draw[->] (0.9,-0.1) to  (0.1,-0.9);
\draw[->] (-0.9,0.1) to (-0.1,0.9);
\draw[->] (-0.9,-0.1) to  (-0.1,-0.9);

\node[circle,draw,inner sep=1pt,label=right:{$A$}] (A) at (1,0) {};
\node[circle,draw,inner sep=1pt,label=above:{$B$}] (B) at (0,1) {};
\node[circle,draw,inner sep=1pt,label=left:{$A'$}] (A') at (-1,0) {};
\node[circle,draw,inner sep=1pt,label=below:{$B'$}] (B') at (0,-1) {};

\end{tikzpicture}
\caption{Any transport plan between $\frac12 (\delta_A+\delta_{A'})$ and $\frac12(\delta_B+\delta_{B'})$ is optimal} \label{fig:carre}
\end{figure}
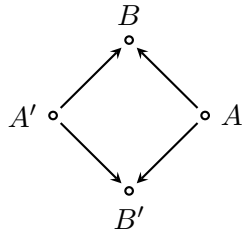
\end{center}

It is not difficult to see that all transport plans have the same cost, in particular there exist infinitely many optimal transport plans. Moreover, if we keep $A$ and $A'$ unchanged, but we move a little bit $B$ and $B'$, then we can recover unique optimal transport plans/maps. For instance, if we slightly move $B$ horizontally to the right, and $B'$ symmetrically slightly to the left, then we recover a unique optimal transport map, where the mass at $A$ is sent to $B$, and the mass at $A'$ is sent to $B'$. Symmetrically, if we move $B$ to the left and  $B'$ to the right, then we also recover a unique optimal transport map, but the mass at $A$ is sent to $B'$ and the mass at $A'$ is sent to $B$. The configuration displayed in Figure \ref{fig:carre} therefore generates a strong instability of optimal transport maps (and plans). However, this example is quite specific in the sense that it deals with discrete measures, and Brenier's theorem does not apply. 

Nonetheless, we can take inspiration from this example to construct absolutely continuous (and compactly supported) source measures $\rho$ for which optimal transport maps are highly unstable. First, we smoothen a little bit the source measure by replacing it by $\rho_r=\frac{1}{2\pi r^2}(\delta_{B(A,r)}+\delta_{B(A',r)})$ for some small $r>0$. We denote by $T^{(r)}_\mu$ the optimal transport map from $\rho_r$ to a measure $\mu$. We know by the results of \cite{letrouitmerigot} that stability holds for any $r>0$: for any compact set $\Y\subset\R^d$,
$$
\forall \mu,\nu\in\mathcal{P}(\Y),\qquad \|T^{(r)}_\mu - T^{(r)}_\nu\|_{L^2(\rho)} \leq C_r W_1(\mu,\nu)^{1/6}.
$$
It can be checked that the constant $C_r$ tends to $+\infty$ as $r\rightarrow 0$, which is a sign of instability (in the limit $r\rightarrow 0$). To prove Theorem \ref{t:highlyunstable} we repeat the above construction ``at all scales". In other words, we take $\rho$ uniform over an infinite union of well-chosen pairs of balls $B(A_i,r_i)$, $B(A_i',r_i)$ of various sizes. Actually, we do not take balls but rectangular parallelepipeds because it is helpful at some point to have more parameters than just the radius.

\subsection{The construction} 
%We start by constructing a probability measure $\rho$ which satisfies all conclusions of Theorem \ref{t:unstable}, except that its support is the closure of a bounded open which is \emph{not pathwise connected}.

We first construct the support of $\rho$, denoted by $\X$ in what follows (and depicted in Figure \ref{fig:supp}). In the sequel,
we consider the rectangular parallelepiped
$$
\mathcal{Q}(\ell, r)=\left\{(x_1,\ldots,x_d)\in\R^d \mid 0< x_1<\ell,\ |x_2|< \frac{r}{2}, \text{ and }  |x_i|< \frac12 \text{ for $3\leq i\leq d$}  \right\}
$$
and consider for $A\in\R^d$ the translates
$$
\mathcal{T}^+(A,\ell, r)=A+\mathcal{Q}(\ell,r), \qquad \mathcal{T}^-(A,\ell,r)=A-\mathcal{Q}(\ell,r)
$$
(in the sense of the Minkowski sum and difference).

Let $(\ell_i),(r_i),(w_i)\in(\R_+)^\N$ and $(u_i)\in\R^\N$. We will make several assumptions on these sequences in Sections \ref{s:prelim} and \ref{s:unstabW1}. Let
$$
\mathcal{S}_i=\mathcal{T}^+(A_i^+,\ell_i,r_i)\cup \mathcal{T}^-(A_i^-,\ell_i,r_i)
$$
where $A_i^+=(u_i+w_i,0,\ldots,0)$ and $A_i^-=(u_i-w_i,0,\ldots,0)$, and let
$$
\X= \bigcup_{i=1}^{+\infty} \mathcal{S}_i.
$$
In the sequel, the sequences will always be chosen in a way that when traveling along the $x_1$-axis in the increasing $x_1$ direction, the parallelepipeds defined above do not meet, and are in the order $\mathcal{S}_1,\mathcal{S}_2,\ldots$ as on Figure \ref{fig:supp} below.
%Since $v_i>r_i$, $\X$ is the union of a set of ball contained in the upper half-plane $\{x_2>0\}$, with its symmetrization with respect to the horizontal axis $\{x_2=0\}$. 

\begin{center}
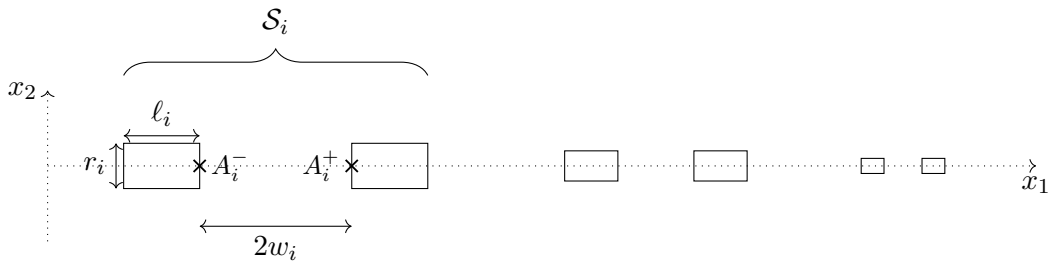
\begin{figure}[h]
\begin{tikzpicture}

  \tikzset{
    axis/.style = {line width=0.8pt},
    tick/.style = {line width=1pt},
    bigcircle/.style = {line width=0.9pt},
    smallsym/.style = {line width=0.7pt},
    lab/.style = {font=\small}
  }

  % La ligne horizontale
  %\draw (-1,0) -- (8,0);
\draw[dotted,->] (-2,-1) -- (-2,1);

  % Nom de l'axe
  \node[left] at (-2,1) {$x_2$};

\draw[dotted,->] (-2,0) -- (11,0);

  % Nom de l'axe
  \node[below] at (11,0) {$x_1$};

  % Les cercles (rayons décroissants)
%  \draw (0,0) circle (0.5);   % grand cercle
 % \draw (2,0) circle (0.5);
  %\draw (5.5,0) circle (0.2);
%  \draw (6.5,0) circle (0.2);
 % \draw (9,0) circle (0.1);
  %\draw (9.5,0) circle (0.1);

%\draw[->] (0,0) -- (0.35,0.35) node[midway,above] {$r_i$};

%  \draw (0,0) circle (0.5);   % grand cercle
 \draw (0,0.3)--(0,-0.3)--(-1,-0.3)--(-1,0.3)--cycle;
 \draw (2,0.3)--(2,-0.3)--(3,-0.3)--(3,0.3)--cycle;

 \draw (5.5,0.2)--(5.5,-0.2)--(4.8,-0.2)--(4.8,0.2)--cycle;
 \draw (6.5,0.2)--(6.5,-0.2)--(7.2,-0.2)--(7.2,0.2)--cycle;

 \draw (9,0.1)--(9,-0.1)--(8.7,-0.1)--(8.7,0.1)--cycle;
 \draw (9.5,0.1)--(9.5,-0.1)--(9.8,-0.1)--(9.8,0.1)--cycle;

\draw[<->] (-1.1,-0.3) -- (-1.1,0.3) node[midway,left] {$r_i$};
\draw[<->] (-1,0.4) -- (0,0.4) node[midway,above] {$\ell_i$};
\draw[<->] (0,-0.8) -- (2,-0.8) node[midway,below] {$2w_i$};

  \coordinate (A1) at (0,0);
  \draw[smallsym] ($(A1)+(-0.08,0.08)$) -- ($(A1)+(0.08,-0.08)$);
  \draw[smallsym] ($(A1)+(-0.08,-0.08)$) -- ($(A1)+(0.08,0.08)$);
  \node[lab, right=0.5pt] at ($(A1)$) {$A_i^-$};

  \coordinate (A1p) at (2,0);
  \draw[smallsym] ($(A1p)+(-0.08,0.08)$) -- ($(A1p)+(0.08,-0.08)$);
  \draw[smallsym] ($(A1p)+(-0.08,-0.08)$) -- ($(A1p)+(0.08,0.08)$);
  \node[lab, left=0.5pt] at ($(A1p)$) {$A_i^+$};

\draw [decorate,decoration={brace,amplitude=10pt}]
    (-1,1.2) -- (3,1.2)
    node[midway,above=12pt] {$\mathcal{S}_i$};

\end{tikzpicture}
\caption{Part of the support of $\rho$, projected on the $(x_1,x_2)$-plane} \label{fig:supp}
\end{figure}
\end{center}

We consider $\rho$ an absolutely continuous probability measure whose support is $\X$ and which, for any $i$, is uniform on $\mathcal{S}_i$. Finally, we let for any $i$
$$
\sigma_i=\rho(\mathcal{T}^+(A_i^+,\ell_i,r_i))=\rho(\mathcal{T}^-(A_i^-,\ell_i,r_i)).
$$ %Of course its density is bounded above and below on $\X$ by positive constants. 
With a slight abuse of notation, the density of $\rho$ with respect to Lebesgue is also denoted by $\rho$ in the sequel.

Let $B_i^+=(u_i,w_i,0,\ldots,0)$ and $B_i^-=(u_i,-w_i,0,\ldots,0)$. Let
\begin{equation}\label{e:defmu}
\mu=\sum_{i=1}^{+\infty} \sigma_i (\delta_{B_i^+}+\delta_{B_i^-}).
\end{equation}
It is immediate to check that this is a probability measure.

For $i\in\N$, let $C_i^+=(u_i+r_i,w_i,0,\ldots,0)$, $C_i^-=(u_i-r_i,-w_i,0,\ldots,0)$, and
$$
\nu_i=\mu+\sigma_i(\delta_{C_i^+}+\delta_{C_i^-}-\delta_{B_i^+}-\delta_{B_i^-}).
$$
It is of course a probability measure.
The only difference between $\nu_i$ and $\mu$ is that in the sum \eqref{e:defmu}, the $i$-th term has been replaced by $\sigma_i(\delta_{C_i^+}+\delta_{C_i^-})$, while all other terms are left unchanged. With the choices made in the next section, $C_i^+$ and $C_i^-$ are seen as perturbations of $B_i^+$ and $B_i^-$.

In the sequel we call ``the $i$-th cell" the set 
$$
\pi_{12}(\mathcal{S}_i\cup \{B_i^+, B_i^-, C_i^+, C_i^-\})
$$
where $\pi_{12}$ denotes the projection on the first two coordinates.

\subsection{Preliminary computations.} \label{s:prelim}
In the sequel we take the convention $u_0=-\infty$.
We make the following assumptions on the sequences (roughly illustrated in Figure \ref{fig:supp2}): for any $i\in\N$, 
\begin{align}
\min(u_i-u_{i-1},u_{i+1}-u_i)& \geq 100 \max(\ell_i,r_i,w_i) \label{e:tailleeps}\\
%\varepsilon_i&=\frac{w_ir_i}{v_i}\label{e:vraitailleeps}\\
w_i&\geq 100r_i. \label{e:viri}
\end{align}

The first inequality means that the cells are well-separated: the distances between points inside a given cell are much smaller than the distance to the closest other cells. The second inequality means $B_i^+$ is much closer to $C_i^+$ than to $B_i^-$ (and $B_i^-$ is much closer to $C_i^-$ than to $B_i^+$). We will choose specific sequences $(\ell_i)$, $(r_i)$, $(w_i)$ and $(u_i)$ in Section \ref{s:unstabW1}. In any case of application, these sequences are bounded (which implies that the probability measure $\rho$ which we construct has compact support).

\begin{center}
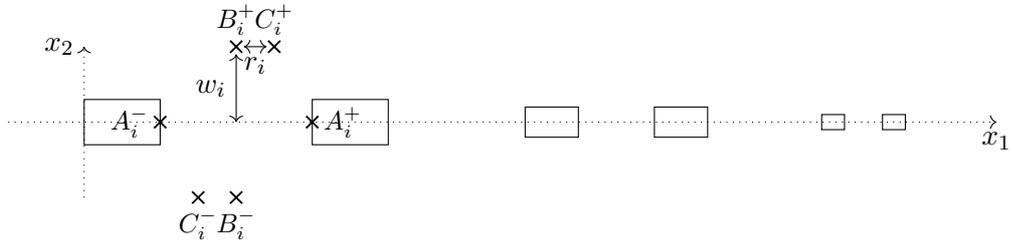
\begin{figure}[h]
\begin{tikzpicture}

  \tikzset{
    axis/.style = {line width=0.8pt},
    tick/.style = {line width=1pt},
    bigcircle/.style = {line width=0.9pt},
    smallsym/.style = {line width=0.7pt},
    lab/.style = {font=\small}
  }

  % La ligne horizontale
  %\draw (-1,0) -- (8,0);
\draw[dotted,->] (-1,-1) -- (-1,1);

  % Nom de l'axe
  \node[left] at (-1,1) {$x_2$};

\draw[dotted,->] (-2,0) -- (11,0);

  % Nom de l'axe
  \node[below] at (11,0) {$x_1$};

  % Les cercles (rayons décroissants)
%  \draw (0,0) circle (0.5);   % grand cercle
 \draw (0,0.3)--(0,-0.3)--(-1,-0.3)--(-1,0.3)--cycle;
 \draw (2,0.3)--(2,-0.3)--(3,-0.3)--(3,0.3)--cycle;

 \draw (5.5,0.2)--(5.5,-0.2)--(4.8,-0.2)--(4.8,0.2)--cycle;
 \draw (6.5,0.2)--(6.5,-0.2)--(7.2,-0.2)--(7.2,0.2)--cycle;

 \draw (9,0.1)--(9,-0.1)--(8.7,-0.1)--(8.7,0.1)--cycle;
 \draw (9.5,0.1)--(9.5,-0.1)--(9.8,-0.1)--(9.8,0.1)--cycle;

\draw[<->] (1,0) -- (1,0.9) node[midway,left] {$w_i$};
\draw[<->] (1.1,1) -- (1.4,1) node[midway,below] {$r_i$};
%\draw[<->] (0,-0.8) -- (2,-0.8) node[midway,below] {$2v_i$};

  \coordinate (A1) at (0,0);
  \draw[smallsym] ($(A1)+(-0.08,0.08)$) -- ($(A1)+(0.08,-0.08)$);
  \draw[smallsym] ($(A1)+(-0.08,-0.08)$) -- ($(A1)+(0.08,0.08)$);
  \node[lab, left=0.5pt] at ($(A1)$) {$A_i^-$};

  \coordinate (A1p) at (2,0);
  \draw[smallsym] ($(A1p)+(-0.08,0.08)$) -- ($(A1p)+(0.08,-0.08)$);
  \draw[smallsym] ($(A1p)+(-0.08,-0.08)$) -- ($(A1p)+(0.08,0.08)$);
  \node[lab, right=0.5pt] at ($(A1p)$) {$A_i^+$};

  \coordinate (B1) at (1,1);
  \draw[smallsym] ($(B1)+(-0.08,0.08)$) -- ($(B1)+(0.08,-0.08)$);
  \draw[smallsym] ($(B1)+(-0.08,-0.08)$) -- ($(B1)+(0.08,0.08)$);
  \node[lab, above=0.5pt] at ($(B1)$) {$B_i^+$};

  \coordinate (B1p) at (1,-1);
  \draw[smallsym] ($(B1p)+(-0.08,0.08)$) -- ($(B1p)+(0.08,-0.08)$);
  \draw[smallsym] ($(B1p)+(-0.08,-0.08)$) -- ($(B1p)+(0.08,0.08)$);
  \node[lab, below=0.5pt] at ($(B1p)$) {$B_i^-$};

  \coordinate (C1) at (1.5,1);
  \draw[smallsym] ($(C1)+(-0.08,0.08)$) -- ($(C1)+(0.08,-0.08)$);
  \draw[smallsym] ($(C1)+(-0.08,-0.08)$) -- ($(C1)+(0.08,0.08)$);
  \node[lab, above=0.5pt] at ($(C1)$) {$C_i^+$};

  \coordinate (C1p) at (0.5,-1);
  \draw[smallsym] ($(C1p)+(-0.08,0.08)$) -- ($(C1p)+(0.08,-0.08)$);
  \draw[smallsym] ($(C1p)+(-0.08,-0.08)$) -- ($(C1p)+(0.08,0.08)$);
  \node[lab, below=0.5pt] at ($(C1p)$) {$C_i^-$};
\end{tikzpicture}
\caption{The support of the measures $\rho$, $\mu$ and $\nu_i$} \label{fig:supp2}
\end{figure}
\end{center}

Under assumptions \eqref{e:tailleeps}-\eqref{e:viri}, we show that for any choice of $C,\alpha>0$ and $p\geq 1$, the inequality 
\begin{equation}\label{e:stabmunui}
\|T_\mu-T_{\nu_i}\|_{L^2(\rho)} \leq CW_p(\mu,\nu_i)^\alpha
\end{equation}
cannot hold for all $\nu_i$ simultaneously. For this, we first compute $T_\mu$ and $T_{\nu_i}$ for any $i$.

For $x=(x_1,\ldots,x_d)\in\R^d$, let us check that
\begin{equation}\label{e:fauxTmu}
x\mapsto \begin{cases} 
B_i^+ & \text{if\ }  x\in \mathcal{S}_i  \text{ \ and \ } x_2\geq 0 \\   B_i^- &  \text{if\ }  x\in \mathcal{S}_i  \text{ \ and \ } x_2<0 \end{cases}
\end{equation}
coincides $\rho$-a.e. with $T_\mu$.
It is immediate to verify that this application transports $\rho$ to $\mu$. Let us check that each point of the support of $\rho$ is sent by this application to the closest point in the support of $\mu$. Let $x\in \mathcal{S}_i$. We may assume that $x_2\geq 0$ (the case $x_2<0$ is symmetric) and that $x_3=\ldots=x_d=0$ because all points in the support of $\mu$ have these coordinates equal to $0$. Then $x$ is closer to $B_i^+$ than to $B_i^-$, and $x$ is at distance at most $(w_i^2+(\ell_i+w_i)^2)^{1/2}$ from $B_i^+$. Let $y$ be another point in the support of $\mu$, different from $B_i^-$ and $B_i^+$. The triangle inequality and then \eqref{e:tailleeps} yield
\begin{equation}\label{e:xy}
|x-y| \geq \min(u_i-u_{i-1},u_{i+1}-u_i)-w_i-\ell_i > (w_i^2+(\ell_i+w_i)^2)^{1/2} \geq |x-B_i^+|
\end{equation}
i.e., $x$ is closer to $B_i^+$ than to $y$.
Therefore, the transport is made at smallest possible cost, and \eqref{e:fauxTmu} is the optimal transport map.

Regarding $T_{\nu_i}$, let us show that it coincides $\rho$-a.e. with the application
\begin{equation}\label{e:ersatz}
x\mapsto 
\begin{cases} T_\mu(x) & \text{if \ } x\notin \mathcal{S}_i \\ C_i^+ & \text{if \ } x\in \mathcal{T}^+(A_i^+,\ell_i,r_i)  \\ C_i^- & \text{if \ } x\in \mathcal{T}^-(A_i^-,\ell_i,r_i).  \end{cases}
\end{equation}
It is clear that this application defines a transport map from $\rho$ to $\nu_i$. We only need to show the following claim, which implies that \eqref{e:ersatz} is the optimal transport map from $\rho$ to $\nu_i$:

\emph{Claim:} The application \eqref{e:ersatz} sends each point $x$ in the support of $\rho$ to its closest point in the support of $\nu_i$. 

This claim is straightforward to check for $x\in \mathcal{S}_j$ ($j\neq i$) with a similar argument as in \eqref{e:xy}. It is also immediate that the closest point to $x\in \mathcal{T}^+(A_i^+,\ell_i,r_i)$ is either $C_i^+$ or $C_i^-$. Let us show that any $x\in \mathcal{T}^+(A_i^+,\ell_i,r_i)$ is closer to $C_i^+$ than to $C_i^-$. We let $x=A_i^++(x_1,\ldots,x_d)$ and observe that
\begin{align*}
|x-C_i^-|^2 -|x-C_i^+|^2&=  (x_1+w_i+r_i)^2 + (w_i+x_2)^2 - (x_1+w_i-r_i)^2 - (w_i-x_2)^2 \\
&=4r_i (w_i+x_1)+4w_ix_2\\
&\geq 4r_iw_i-2r_iw_i\\
&>0.
\end{align*}
Similarly if $x\in \mathcal{T}^-(A_i^-,\ell_i,r_i)$, then $x$ is closer to $C_i^-$ than to $C_i^+$. Hence the transport is made at smallest possible cost, and \eqref{e:ersatz} is the optimal transport map.

%\textcolor{red}{Dire que chaque point est envoyé vers le point le plus proche.} Let $\pi_1$ denote the projection onto the first coordinate, and consider an auxiliary measure $\rho'$ which is the pushforward of $\rho$ under $\pi_1$. The optimal transport from $\rho'$ to $\mu$ is a 1d problem, so there is a convex function $\widetilde{\phi}$ on $\R$ such that $\nabla\widetilde{\phi}$ is the optimal transport map from $\rho'$ to $\mu$. We let $\phi(x_1,x_2)=\widetilde{\phi}(x_1)$, which is also a convex function, and we observe that $(\nabla\phi)_{\#}\rho = \mu$. Therefore $T_\mu=\nabla\phi$ is the optimal transport map from $\rho$ to $\mu$. We check immediately that the points in $B(A_i,r_i)$ which are located on the left of $A_i$ are sent to $B_i$ and the points on its right are sent to $B_i'$.

%Let us now describe the optimal transport map from $\rho$ to $\nu_i$. Forgetting about the part of the measures around figure $i$, we need to have an optimal transport, therefore the optimal transport from $\rho$ to $\nu_i$ is as described above. For the part around the $i$-th figure, since any point in $B(A_i,r_i)$ is closer to $C_i$ and any point in $B(A_i',r_i)$ is closer to $C_i'$, the optimal transport map is to send $B(A_i,r_i)$ to $C_i$ and $B(A_i',r_i)$ to $C_i'$.

From the explicit expressions of $T_\mu$ and $T_{\nu_i}$ obtained above, we deduce
\begin{equation}\label{e:exactTmuTnu}
\|T_\mu-T_{\nu_i}\|^2_{L^2(\rho)} =\frac12 \rho(\mathcal{S}_i)|B_i^+-C_i^+|^2 +\frac12 \rho(\mathcal{S}_i)|B_i^+-C_i^-|^2 = \sigma_i(2r_i^2+4w_i^2)\geq 4\sigma_iw_i^2.
\end{equation}
Besides, due to \eqref{e:viri}, we have $r_i\leq w_i/100$, hence the optimal coupling between $\mu$ and $\nu_i$ is a coupling where the mass at $B_j^+$ (resp. $B_j^-$) does not move if $j\neq i$, and is sent to $C_j^+$ (resp. $C_j^-$) if $j=i$. Indeed, in this way, each piece of mass in the support of $\mu$ is sent to the closest point in the support of $\nu_i$. Therefore, for any  $p\geq 1$,
\begin{equation}\label{e:exactWp}
W_p(\mu,\nu_i) = r_i \sigma_i^{1/p}.
\end{equation}
Combining \eqref{e:exactTmuTnu} and \eqref{e:exactWp} we deduce that for any $\alpha>0$, %Recalling that $\sigma_i=\omega_d r_i^d$ where $\omega_d$ is the volume of the $d$-dimensional unit ball, 
\begin{equation}\label{e:impequiv}
\frac{\|T_\mu-T_{\nu_i}\|^2_{L^2(\rho)}}{W_p(\mu,\nu_i)^{2\alpha}}\geq 4 w_i^{2}r_i^{-2\alpha}\sigma_i^{1-\frac{2\alpha}{p}}.
\end{equation}
%where the inequality comes from $\alpha\geq 1$ and $r_i<1$, and the convergence is due to our assumption \eqref{e:grandedistver}. This shows that \eqref{e:stabmunui} does not hold, for any $C>0$.

The end of the paper is devoted to the conclusion of the proofs of Theorems \ref{t:highlyunstable} and \ref{t:unstableW1}.

%\subsection{Proof of Theorem \ref{t:noncompacttarget}} \label{s:noncompacttarget}We choose $r_i$ decaying exponentially fast, $\ell_i=r_i$, and $\sigma_i=r_i^2$. We also set $w_i=i^{-10}r_i^{-1}$ and let $v_i=i^{-2}$. \textcolor{red}{A reprendre, j'ai changé des exposants pour que $\varepsilon_i$ soit vraiment petit, là c'est $i^{-8}$.} Then $\mu\in\mathcal{P}_2(\R^d)$since$$\int_{\R^d}|x|^2\dd\mu \leq \sum_i \sigma_i w_i^2 = \sum_i i^{-20} <+\infty.$$Moreover,\begin{equation}\label{e:impequiv3}\frac{\|T_\mu-T_{\nu_i}\|^2_{L^2(\rho)}}{W_p(\mu,\nu_i)^{2\alpha}}\geq \frac{4w_i^2\sigma_i}{(\frac{2w_ir_i}{v_i}\sigma_i^{1/p})^{2\alpha}}= \frac{4i^{-2}}{(2i^{3}r_i^{2/p})^{2\alpha}}.\end{equation}Thus we need to prove that for any $n\in\N$, $$r_i^{2/p}\ll i^{-n}.$$This holds due to our choices. Notice that $W_p(\mu,\nu_i)=\frac{2w_ir_i}{v_i}\sigma_i^{1/p}=2i^3 r_i^{2/p}$ tends to $0$ as $i\rightarrow +\infty$. 

\subsection{End of the proof of Theorem \ref{t:unstableW1}} \label{s:unstabW1}
Let $p\geq 1$ and $\alpha>\frac{p}{2(p+1)}$. We choose the sequences $(\ell_i)_{i\in\N}$, $(r_i)_{i\in\N}$, $(w_i)_{i\in\N}$ and $(u_i)_{i\in\N}$ as follows. We let $r_i=c_02^{-i}$ (again, any sequence with superpolynomial decay would work), $\ell_i=w_i=c_0c_1i^{-2}$ and $u_1\in\R$, $u_{i+1}-u_i=c_0c_2i^{-2}$ where $c_1,c_2>0$ are chosen in a way that \eqref{e:tailleeps} and \eqref{e:viri} hold.
Choosing appropriately $c_0>0$, we fix
$$
\sum_{i=1}^{+\infty} \ell_i r_i=\frac12
$$
and let $\sigma_i=\ell_ir_i$, which makes $\rho$ uniform over its support. We get
\begin{equation}\label{e:impequiv2}
\frac{\|T_\mu-T_{\nu_i}\|^2_{L^2(\rho)}}{W_p(\mu,\nu_i)^{2\alpha}}\geq 4w_i^2r_i^{-2\alpha}\sigma_i^{1-\frac{2\alpha}{p}}=Cw_i^{3-\frac{2\alpha}{p}}r_i^{1-\frac{2\alpha}{p}-2\alpha}.
\end{equation}
Since $\alpha>\frac{p}{2(p+1)}$, this quantity tends to $+\infty$ as $i\rightarrow +\infty$.
%We see that for $\alpha>\frac{d}{2(d+1)}$, the exponent of $r_i$ is negative, and therefore choosing appropriate sequences satisfying \eqref{e:tailleeps} we get that the quantity \eqref{e:equivforunif} which tends to $+\infty$ as $i\rightarrow +\infty$.

\begin{remark}
We did not obtain any improvement of our results by introducing other parameters in the geometric picture, for instance by making the distance between $A_i^-$ and $A_i^+$ different from the distance between $B_i^-$ and $B_i^+$.
\end{remark}

%\begin{remark}We could take the second coordinate of $B_i$ to be different from the half-distance between $A_i$ and $A_i'$, and optimize over this additional parameter, but we did not manage to improve the result using this additional parameter.\end{remark}

\end{document}